# Assessing Nonlinear Diffusion Acceleration for Boltzmann Fokker Planck Equation in Slab Geometry

Japan K. Patel[1,*], Barry D. Ganapol[2], Martha M. Matuszak[1]

[1]University of Michigan, Ann Arbor, MI; [2]University of Arizona, Tucson, AZ

## ABSTRACT

The convergence of Boltzmann Fokker Planck solution can become arbitrarily slow with iterative procedures like source iteration. This paper derives and investigates nonlinear diffusion acceleration scheme for the solution of Boltzmann Fokker Planck equation in slab geometry. This method is a conventional high order/low order scheme with a traditional "diffusion-plus-drift" low order system. The method, however, differs from the earlier variants as the definition of the low order equation, which is adjusted according to the zeroth and first moments of the Boltzmann Fokker Planck equation. For the problems considered, we observe that the NDA-accelerated solution follows the unaccelerated well and provides roughly an order of magnitude savings in iteration count and runtime compared to source iteration.

## 1. INTRODUCTION

Transport problems with highly forward-peaked scattering are encountered in several applications including radiation shielding, medical physics, and plasma physics. Such problems often have small mean free paths and nearly singular differential scattering cross-sections in the forward direction, and therefore require large expansion orders to adequately represent scattering [1]. Several approximation techniques have been proposed to overcome this issue. For example, Fokker-Planck [2], generalized Fokker Planck [3,4], Boltzmann Fokker Planck (BFP) [1,5], and generalized Boltzmann Fokker Planck [6] methods have been investigated. This paper primarily considers the Boltzmann Fokker Planck approximation [1]. This equation has been applied to applications related to both neutral and charged particles [1,5] and treats scattering by decomposing it into its smooth and singular components [1,5]. While the smooth component is represented using Legendre moments, the singular component is represented using the Fokker-Planck operator [1].

Just like the Boltzmann equation, convergence of the solution of the BFP can become arbitrarily slow depending on the cross-sections. In order to speed up convergence, diffusion-based acceleration techniques like diffusion synthetic acceleration (DSA) [6] and quasi-diffusion (QD) [7] are often employed. Nonlinear high-order low-order (HOLO) acceleration techniques provide substantial performance gains for several transport-related single and multiphysics problems [8,9]. These algorithms accelerate the convergence of slowly converging physics, like scattering, to reduce the number of iterations required. They also allow us to isolate high-order (HO) equations from the coupled system (multiphysics or otherwise) by introducing a coarser scale low-order (LO) system. These methods usually cater to isotropic or weakly anisotropic problems, with a few exceptions [7].

---
[*]Corresponding author, jakpatel@umich.edu, Tel: 614-401-0603

In this paper, we implement and investigate the application of the nonlinear diffusion acceleration (NDA) [8] for the solution of the Boltzmann Fokker-Planck equation in slab geometry. We briefly review NDA in the next section and conclude in Section 3.

## 2. NONLINEAR DIFFUSION ACCELERATION

Forward-peaked transport problems are often represented using the BFP approximation, where the scattering kernel is divided into smooth and singular components [1]. With angular flux $\psi$, $l^{th}$ Legendre moment of angular flux $\phi_l$, total macroscopic cross-section $\sigma_t$, $l^{th}$ order Legendre polynomial $P_l$, scattering cross-section moments $\sigma_{s,l}$, and an internal source $Q$, we consider the following BFP equation [1]:

$$\mu \frac{d}{dx}\psi(x,\mu) + \sigma_t(x)\psi(x,\mu) = \sum_{l=0}^{L-2} \frac{(2l+1)}{2} P_l(\mu)\tilde{\sigma}_{s,l}(x)\phi_l(x) + \frac{\sigma_{tr}}{2}\frac{d}{d\mu}(1-\mu^2)\frac{d}{d\mu}\psi + \frac{Q(x,\mu)}{2}. \quad (1)$$

The smooth moments $\tilde{\sigma}_{s,l}$ are evaluated using the original scattering cross-section moments and expansion order $L$ according to [1]:

$$\tilde{\sigma}_{s,l} = \sigma_{s,l} - \sigma_L - \frac{\alpha}{2}\big(L(L+1) - l(l+1)\big), \quad (2)$$

The transfer cross-section $\sigma_{tr}$ is defined as [1]:

$$\sigma_{tr} = \frac{\sigma_{s,L-1} - \sigma_{s,L}}{L}, \quad (3)$$

and moments of the angular flux are determined using the following relation [10]:

$$\phi_l(x) = \int_{1}^{1} \psi(x,\mu) P_l(\mu) d\mu. \quad (4)$$

The paper primarily considers monoenergetic problems, but this approach can be extended to multi-group settings.

### 2.1. Derivation

While this method is well-established and has been investigated at length in many publications [8,9], we present a brief derivation of NDA for the BFP equation. This preliminary work only considers monoenergetic problems and restricts attention to isotropic internal sources. Following [8], and dropping the notation for angle and space dependence, we begin by developing the first two moment equations. The $0^{th}$ moment equation is obtained by multiplying the BFP equation by $P_0$ and integrating over the angular space:

$$\int_{-1}^{1} d\mu \, P_0 \left( \mu \frac{d}{dx}\psi + \sigma_t \psi = \sum_{l=0}^{L-2} \frac{(2l+1)}{2} P_l \tilde{\sigma}_{s,l} \phi_l + \frac{\sigma_{tr}}{2}\frac{d}{d\mu}(1-\mu^2)\frac{d}{d\mu}\psi + \frac{Q}{2} \right). \quad (5)$$

Using the recurrence relation for Legendre polynomials [10], and the definition of moment of $0^{th}$ and $1^{st}$ moments from Eq. 4, Eq. 5 can be simplified to:

$$\frac{d}{dx}\phi_1 + \sigma_t\phi_0 = \tilde{\sigma}_{s,0}\phi_0 + \int_{-1}^{1} d\mu \, P_0 \left(\frac{\sigma_{tr}}{2}\frac{d}{d\mu}(1-\mu^2)\frac{d}{d\mu}\psi\right) + Q_0. \tag{6}$$

We now expand angular flux in the fourth term of Eq. 6 using Legendre polynomials [10], rearrange the equation, and rewrite the angular Laplacian (Fokker-Planck) operator according to its eigenvalues and eigenfunctions [11] to get:

$$\frac{d}{dx}\phi_1 + \sigma_t\phi_0 = \tilde{\sigma}_{s,0}\phi_0 + \int_{-1}^{1} d\mu \, P_0 \sum_{l}^{l} \frac{2l+1}{2}\phi_l \frac{\sigma_{tr}}{2}(-l(l+1))P_l + Q_0. \tag{7}$$

Rearranging the fourth term of Eq. 7 and using the orthogonality property of Legendre polynomials returns the following $0^{th}$ moment equation:

$$\frac{d}{dx}\phi_1 + (\sigma_t - \tilde{\sigma}_{s,0})\phi_0 = Q_0. \tag{8}$$

Next, we multiply the BFP equation with $P_1$ and integrate it over angular space:

$$\int_{-1}^{1} d\mu \, P_1 \left(\mu\frac{d}{dx}\psi + \sigma_t\psi = \sum_{l=0}^{L-2}\frac{(2l+1)}{2}P_l\tilde{\sigma}_{s,l}\phi_l + \frac{\sigma_{tr}}{2}\frac{d}{d\mu}(1-\mu^2)\frac{d}{d\mu}\psi + \frac{Q}{2}\right). \tag{9}$$

Using the same steps as the ones used to derive the zeroth moment equation, we obtain the following first moment equation:

$$\frac{1}{3}\frac{d}{dx}\phi_0 + \frac{2}{3}\frac{d}{dx}\phi_2 + (\sigma_t - \tilde{\sigma}_{s,1} + \sigma_{tr})\phi_1 = 0. \tag{10a}$$

We then replace the eliminate the term containing $\phi_2$ and add a drift term to preserve consistency in Eq. 10a to obtain the following relation for the current $\phi_1$:

$$\phi_1 = -\frac{1}{3(\sigma_t - \tilde{\sigma}_{s,1} + \sigma_{tr})}\frac{d}{dx}\phi_0 + \widehat{D}\phi_0. \tag{10b}$$

This diffusion-plus-drift relation for the current is substituted back into Eq. 8 to obtain the following LO equation:

$$\frac{d}{dx}\left(-\frac{1}{3(\sigma_t - \tilde{\sigma}_{s,1} + \sigma_{tr})}\frac{d}{dx}\phi_0 + \widehat{D}\phi_0\right) + (\sigma_t - \sigma_{s,0})\phi_0 = Q, \tag{11}$$

Therefore, HOLO system takes the following form:

$$\text{HO: } \mu\frac{d}{dx}\psi^{HO} + \sigma_t\psi^{HO} = \sum_{l=0}^{L-2}\frac{(2l+1)}{2}P_l\tilde{\sigma}_{s,l}\phi_l + \frac{\sigma_{tr}}{2}\frac{d}{d\mu}(1-\mu^2)\frac{d}{d\mu}\psi + \frac{Q}{2},$$

$$\text{LO: } -\frac{d}{dx}\left(\frac{1}{3(\sigma_t - \tilde{\sigma}_{s,1} + \sigma_{tr})}\frac{d}{dx}\phi_0^{LO} + \widehat{D}\phi_0^{LO}\right) + (\sigma_t - \sigma_{s,0})\phi_0^{LO} = Q, \tag{12a}$$

$$\text{Closure: } \widehat{D} = \frac{J^{HO} + \frac{1}{3(\sigma_t - \tilde{\sigma}_{s,1} + \sigma_{tr})}\frac{d\phi^{HO}}{dx}}{\phi^{HO}}, \tag{12b}$$

(12c)

with boundary conditions chosen to ensure consistency between HO and LO systems [8]:

$$\frac{\phi_1^{LO}}{\phi_0^{LO}} = \frac{\phi_1^{HO}}{\phi_0^{HO}}. \tag{13}$$

This nonlinear system of equations has been solved using various numerical techniques [12]. We use Picard iteration [9] in this paper for simplicity.

## 2.2. Discretization and solution

The BFP equation is discretized in angle using the discrete ordinates ($S_N$) method [10]. With quadrature order $N$ and angular index $n$, the BFP-$S_N$ equations are written as:

$$\mu \frac{d}{dx}\psi_n + \sigma_t \psi_n = \sum_{l=0}^{L-2} \frac{(2l+1)}{2} P_l(\mu_n) \tilde{\sigma}_{s,l} \phi_l + \frac{\sigma_{tr}}{2}\left(\frac{d}{d\mu}(1-\mu^2)\frac{d}{d\mu}\psi\right)_n + \frac{Q_n}{2}. \tag{14}$$

We employ Morel's weighted finite difference technique [13] to discretize the Fokker-Planck term in Eq. 14:

$$\left(\frac{d}{d\mu}(1-\mu^2)\frac{d}{d\mu}\psi\right)_n = \gamma_{n+\frac{1}{2}} \dot{\psi}_{n+\frac{1}{2}} - \gamma_{n-\frac{1}{2}} \dot{\psi}_{n-\frac{1}{2}}, \tag{15a}$$

with

$$\dot{\psi}_{n+\frac{1}{2}} = \frac{\psi_{n+1} - \psi_n}{\mu_{n+1} - \mu_n}, \text{ and } \gamma_{n+\frac{1}{2}} = \gamma_{n-\frac{1}{2}} + \nu\mu_n w_n, \tag{15b}$$

where $\nu$ is a normalization constant [13]. The $S_N$ equations are discretized in space using the standard diamond difference technique [13]. We follow Knoll et al. [8] and discretize the LO system using a finite difference approach:

$$-\frac{1}{3(\sigma_t - \tilde{\sigma}_{s,1} + \sigma_{tr})}\frac{\phi_{i+1} - 2\phi_i + \phi_{i-1}}{\Delta x^2} + \frac{\widehat{D}_{i+\frac{1}{2}}}{2\Delta x}(\phi_{i+1} - \phi_i) - \frac{\widehat{D}_{i-\frac{1}{2}}}{2\Delta x}(\phi_i - \phi_{i-1}) + (\sigma_t - \sigma_{s,0})\phi_{0,i} = Q_i, \tag{16a}$$

where $i$ is the spatial discretization index according to Fig. 1, and the consistency factor:

$$\widehat{D}_{i+\frac{1}{2}} = \frac{\phi_{1,i+\frac{1}{2}}^{HO} + \frac{1}{3(\sigma_t - \tilde{\sigma}_{s,1} + \sigma_{tr})}\frac{\phi_{0,i+1}^{HO} - \phi_{0,i}^{HO}}{\Delta x}}{\phi_{0,i+\frac{1}{2}}^{HO}}, \tag{16b}$$

The discretized HOLO system is solved using a simple Picard iteration framework. With iteration index $k$, we sweep through the HO system to get a rough estimate of the angular flux, scalar flux, and current:

$$\mu\frac{d}{dx}\psi_n^{k+\frac{1}{2}} + \sigma_t\psi_n^{k+\frac{1}{2}} - \frac{\sigma_{tr}}{2}\left(\frac{d}{d\mu}(1-\mu^2)\frac{d}{d\mu}\psi\right)_n^{k+\frac{1}{2}} = \sum_{l=0}^{L-2}\frac{(2l+1)}{2}P_l(\mu_n)\tilde{\sigma}_{s,l}\phi_l^{LO,k} + \frac{Q_n}{2}. \tag{17a}$$

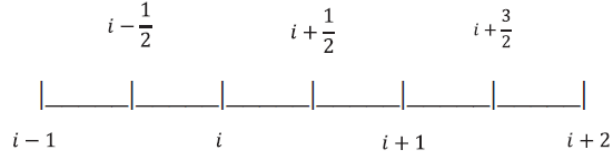

**Figure 1: Spatial grid indexing.**

Subsequently, with total number of nodes $I$ we check for convergence using the following error definition [8]:

$$\epsilon = \frac{1}{\sqrt{I}}\|\phi_0^{HO,k+\frac{1}{2}} - \phi_0^{LO,k}\|, \tag{17b}$$

If necessary, evaluate the consistency factor according to 16b:

$$\widehat{D}_{i+\frac{1}{2}}^{k+1} = \frac{\phi_{1,i+\frac{1}{2}}^{HO,k+\frac{1}{2}} + \frac{1}{3(\sigma_t - \tilde{\sigma}_{s,1} + \sigma_{tr})}\frac{\phi_{0,i+1}^{HO,k+\frac{1}{2}} - \phi_{0,i}^{HO,k+\frac{1}{2}}}{\Delta x}}{\phi_{0,i+\frac{1}{2}}^{HO,k+\frac{1}{2}}}, \tag{17c}$$

and then solve the LO system:

$$-\frac{1}{3(\sigma_t - \tilde{\sigma}_{s,1} + \sigma_{tr})}\frac{\phi_{i+1}^{LO,k+1} - 2\phi_i^{LO,k+1} + \phi_{i-1}^{LO,k+1}}{\Delta x^2} + \frac{\widehat{D}_{i+\frac{1}{2}}^{k+\frac{1}{2}}}{2\Delta x}\left(\phi_{i+1}^{LO,k+1} - \phi_i^{LO,k+1}\right)$$
$$-\frac{\widehat{D}_{i-\frac{1}{2}}^{LO,k+\frac{1}{2}}}{2\Delta x}\left(\phi_i^{LO,k+1} - \phi_{i-1}^{LO,k+1}\right) + (\sigma_t - \sigma_{s,0})\phi_{0,i}^{LO,k+1} = Q_i, \tag{17d}$$

to obtain a revised flux and updated scattering source for the HO sweep. This iteration continues until convergence. We present numerical examples next.

### 2.3. Numerical experiments

In order to obtain a preliminary assessment of the NDA algorithm's performance, we consider two scattering kernels:

1) Henyey-Greenstein kernel (HGK), where the scattering moments are evaluated according to [14]

$$\sigma_{s,l}^{HGK} = \sigma_s g^l. \tag{18}$$

The parameter $g$ is called the anisotropy factor – as it goes from zero to unity, the scattering goes from isotropic to highly anisotropic.

2) Screened Rutherford kernel (SRK), which defines scattering moments as [15]

$$\sigma_{s,l}^{SRK} = \sigma_s \int_{-1}^{1} d\mu P_l \frac{C}{(1+2\eta-\mu)^2}. \tag{19}$$

The parameter $\eta$ defines the degree of anisotropy, with smaller $\eta$ values representing moments with a valid Fokker-Planck limit [15].

We begin our assessment by first comparing source iteration [16] and NDA solutions for the HGK with $g = 0.9$. We set $L = 15$, $N = 16$, and discretize a slab of length 1 cm into 200 spatial cells. We assume a scattering cross-section of unity and a small absorption cross-section of $10^{-6}$ cm$^{-1}$. Additionally, we use a convergence tolerance of $10^{-6}$ and use vacuum boundaries to complete the problem setup. Fig. 2 presents a comparison of scalar fluxes with an increasing number of smooth moments used in the decomposition of the scattering operator. We increase this number from one to thirteen in increments of four. We observe good agreement between unaccelerated and accelerated solutions for all four cases. We note, however, that the accelerated solution does not reproduce the unaccelerated solution to convergence tolerance. This behavior is consistent with nonlinear acceleration techniques like quasi-diffusion (QD) [8], but not with synthetic acceleration techniques like diffusion synthetic acceleration (DSA) [6]. This will be investigated in the future.

We consider one example each from both HGK and SRK to assess the speedup estimates. Table 1 presents a comparison between runtimes and iteration counts observed with both SI and NDA. As expected, we observe roughly one order of magnitude reduction in iteration count and runtime with NDA. This is consistent with speedups DSA provides for highly anisotropic transport problems [11,19]. We note that the iteration count freezes to a fixed value once the number of smooth moments increases beyond 5. This indicates that higher moments do not significantly impact iteration count (although they do contribute to the runtime due to extra computational effort required to obtain them). This, in turn, means that if enough of the lower moments are attenuated, the overall efficiency of the solution algorithm can be improved without having to represent all moments explicitly.

**Table I. Iteration count and runtime comparison.**

| Kernel | Parameter | BFP order | SI iterations | SI runtime [s] | NDA iterations | NDA runtime [s] |
|---|---|---|---|---|---|---|
| HGK | $g = 0.9$ | 1 | 26 | 1.4 | 12 | 0.8 |
| | | 5 | 35 | 4.9 | 18 | 2.6 |
| | | 9 | 36 | 8.2 | 18 | 4.2 |
| | | 13 | 36 | 11.9 | 18 | 5.81 |
| SRK | $C = 0.3903$ $\eta = 2.836 \times 10^{-5}$ | 1 | 2655 | 132.3 | 351 | 18.4 |
| | | 5 | 2739 | 377.4 | 318 | 44.3 |
| | | 9 | 2739 | 626.8 | 318 | 72.2 |

| | | 13 | 2739 | 867.7 | 314 | 100.8 |

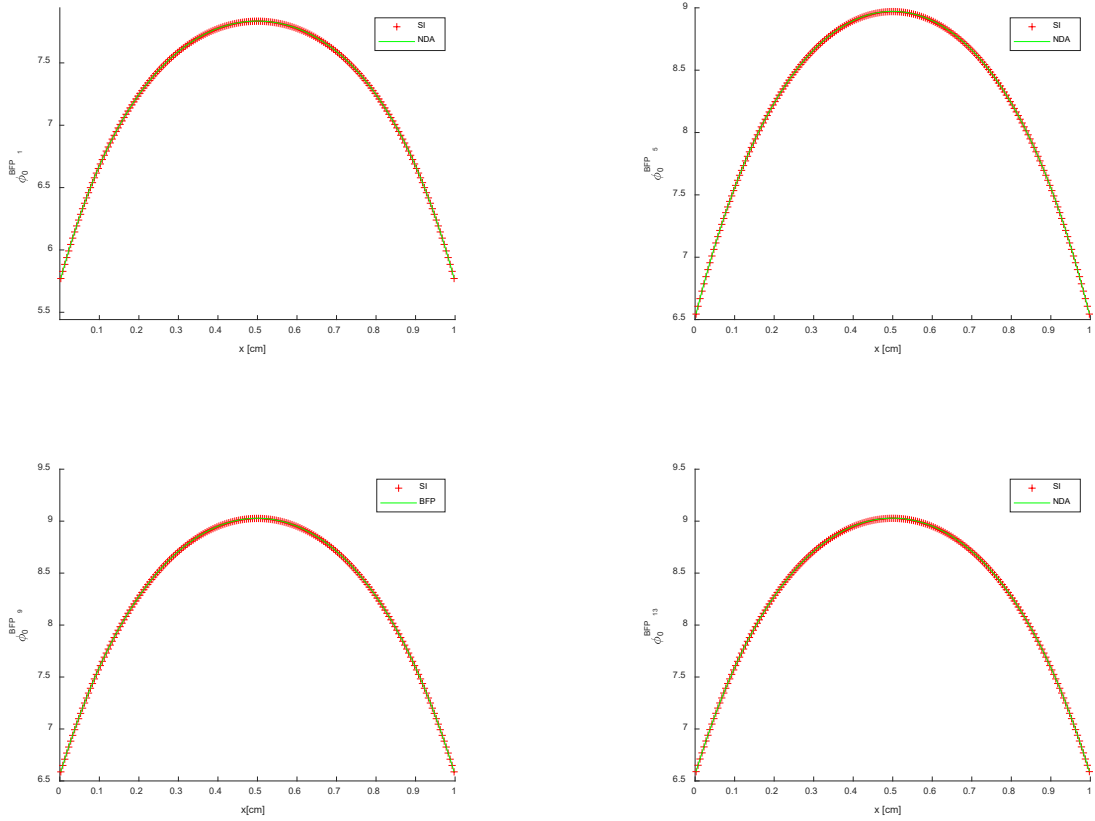

**Figure 2: Comparision of scalar fluxes obtained using SI and NDA with increasing number of smooth moments in BFP decomosition of the scattering operator.**

## 3.  CONCLUSIONS

The BFP approximation decomposes the scattering operator into smooth and singular components, which are represented using a Legendre expansion and an angular Laplacian (Fokker-Planck) operator, respectively [1]. While the NDA scheme is quite well known for the Boltzmann equation [8], it has to be altered to account for the Fokker-Planck term present in the BFP equation. In this abstract, we systematically derived and implemented the nonlinear diffusion acceleration scheme for the BFP equation. We discretized the HO system using the diamond difference method in space and the standard discrete ordinates method in angle [10]. Additionally, we used Morel's weighted difference technique to discretize the Fokker-Planck operator [13]. The LO system employed central finite difference technique for discretization. The overall HOLO system of equations was solved using Picard iteration [12].

We considered two scattering kernels for our preliminary study – HGK and SRK and compared solutions obtained using SI and NDA. For the problems implemented, NDA solutions agreed with those obtained using SI. We observed roughly an order of magnitude speed-up in both overall wall-clock runtime and iteration count. This speedup was expected and was found to be consistent with DSA, as observed in other studies [11]. Our main observation, however, was iteration count stagnates with increasing number of smooth moments used in the BFP approximation. This means that all moments do not need to be accelerated to obtain attractive speedups.

This will be the focus of our future work. A nonlinear $P_l$ acceleration scheme will be developed, and an optimal LO system will be determined. Additionally, a generalized quasi-diffusion approach will also be developed and investigated to represent higher order moments of the BFP equation. While anisotropic sources were not considered in this study, the HOLO scheme will be revised to accommodate isotropic boundary and beam sources.